\documentclass[12pt]{article}

\newcommand{\ol}{\setlength{\itemsep}{0pt.}\begin{enumerate}}
\newcommand{\eol}{\end{enumerate}\setlength{\itemsep}{-\parsep}}

\begin{document}
\bibliographystyle{plain}
\parindent 0cm
\parskip 0.2cm
%\input{psfig}
%%%%%%%%%%%%%%%%%%%%%%%%%%%%%%%%%%%%%%%%%%%%%%%%%%%%%%%%%%%%%%
% log-like functions
\newcommand{\real}{\mbox{real}}
\newcommand{\mod}{\mbox{ mod }}
%%%%%%%%%%%%%%%%%%%%%%%%%%%%%%%%%%%%%%%%%%%%%%%%%%%%%%%%%%%%%%
\newtheorem{theorem}{Theorem}
\newtheorem{THEOREM}{Theorem}[section]
\newtheorem{corollary}{Corollary}
\newtheorem{lemma}{Lemma}
\newtheorem{claim}{Claim}
\newtheorem{exercise}{Exercise}
\newtheorem{definition}{Definition}[section]
\newtheorem{remark}{Remark}[section]
\newtheorem{conjecture}{Conjecture}
\newtheorem{note}{Note}
\newtheorem{example}{Example}[section]
\newtheorem{fact}{Fact}[section]
\newtheorem{PROBLEM}[THEOREM]{Problem}
\newenvironment{problem}
{\begin{PROBLEM} \hspace{-.85em} {\bf :} \rm}
{\end{PROBLEM}}
%%%%%%%%%%%%%%%%%%%%%%%%%%%%%%%%%%%%%%%%%%%%%%%%%%%%%%%%%%%%%%
\font\boldsets=msbm10
\def\G{{\hbox{\boldsets  G}}}
\def\M{{\hbox{\boldsets  M}}}
\def\F{{\hbox{\boldsets  F}}}
\def\N{{\hbox{\boldsets  N}}}
\def\R{{\hbox{\boldsets  R}}}
\def\Z{{\hbox{\boldsets  Z}}}
\def\blackslug{\hbox{\hskip 1pt \vrule width 4pt height 8pt depth 1.5pt
\hskip 1pt}}
\def\emb{embedded }
\def\dope{monotone map }
\def\dopes{monotone maps }
\def\Dope{Monotone Maps }
\def\sgn{{\hbox{sgn}}}

\newcommand{\arrow}[1]{\vec{#1}}
\newcommand {\half}{\frac 1 2}
\newcommand{\ignore}[1]{}

\newcommand{\bbox}{\vrule height7pt width4pt depth1pt}
\def\QED{\quad\blackslug\lower 8.5pt\null\par}
\def\proof{\par\penalty-1000\vskip .5 pt\noindent{\bf Proof\/: }}
%problem
\newcommand{\prb}{\begin{problem}}
\newcommand{\eprb}{\bbox\end{problem}}
%%%%%%%%%%%%%%%%%%%%%%%%%%%%%%%%%%%%%%%%%%%%%%%%%%%%%%%%%%%%%%
\title
{Monotone Maps, Sphericity and Bounded Second Eigenvalue}
\author{
Yonatan Bilu and Nati Linial
\thanks{
Institute of Computer Science,
Hebrew University Jerusalem 91904 Israel
\mbox{\em \{johnblue,nati\}@cs.huji.ac.il}. This research is supported by the Israeli Ministry of Science and the Israel Science Foundation.
}
}
\maketitle

%%%%%%%%%%%%%%%%%%%%%%%%%%%%%%%%%%%%%%%%%%%%%%%%%%%%%%%%%%%%%%

\begin{abstract}
We consider {\em monotone} embeddings of a finite metric space into low dimensional
normed space. That is, embeddings that respect the order among the distances
in the original space. Our main interest is in embeddings into Euclidean spaces.
We observe that any metric on $n$ points can be embedded into $l_2^n$, while,
(in a sense to be made precise later), for almost every $n$-point metric space,
every monotone map must be into a
space of dimension $\Omega(n)$ (Lemma \ref{momad2}).\\
It becomes natural, then, to seek explicit constructions of metric spaces that cannot
be monotonically embedded into spaces of
sublinear dimension. To this end, we employ known results on {\em sphericity} of graphs,
which suggest one example of such a metric space - that defined by a complete bipartite
graph.
We prove that an $\delta n$-regular graph of order $n$, with bounded diameter
has sphericity $\Omega(n/(\lambda_2+1))$, where $\lambda_2$ is the
second largest eigenvalue of the adjacency matrix of the graph, and 
$0 < \delta \leq \half$ is constant
(Theorem \ref{our-bound}).
We also show that while random graphs have linear sphericity, there are
{\em quasi-random} graphs of logarithmic sphericity
(Lemma \ref{alex}).\\
For the above bound to be linear, $\lambda_2$ must be constant.
We show that if the second eigenvalue of an $n/2$-regular graph is bounded by a
constant, then the graph is close to being complete bipartite.
Namely, its adjacency matrix differs from that of a complete bipartite
graph in only $o(n^2)$ entries (Theorem \ref{main1}). Furthermore,
for any $0 < \delta < \half$, and $\lambda_2$, 
there are only finitely many $\delta n$-regular graphs with second eigenvalue
at most $\lambda_2$ (Corollary \ref{no-graphs}).
\end{abstract}

{\em Keywords:} Embedding, Finite Metric Space, Graphs, Sphericity, Eigenvalues, Bipartite Graphs, Second Eigenvalue.
\pagebreak
%%%%%%%%%%%%%%%%%%%%%%%%%%%%%%%%%%%%%%%%%%%%%%%%%%%%%%%%%%%%%%%%%%%

\section{Introduction}
Euclidean embeddings of finite metric spaces have been extensively studied,
with the aim of finding an embedding that doesn't distort the metric too much.
We refer the reader to the survey papers of Indyk (\cite{Indyk}) and
Linial (\cite{Nati}), as well as chapter 15 of Matou{\v{s}}ek's
Discrete Geometry book \cite{Matousek}.
Here we focus on a different type of embeddings. Namely, those that
preserve the order relation of the distances. We call such embeddings {\em monotone}.
There are quite a few applications that make this concept natural and interesting, since
there are numerous algorithmic problems whose solution depends
only on the order among the distances. Specifically, questions that concern
nearest neighbors. The notion of monotone embeddings suggests the following
general strategy toward the resolution of such problems. Namely, embed
the metric space at hand monotonically into a ``nice'' space, for which good
algorithms are known to solve the problem.
Solve the problem in the ``nice'' space - the same solution applies as well for the
original space. ``Nice'' often means a low dimensional normed space.
Thus, we focus on the minimal dimension which permits a monotone embedding.\\
In section \ref{dope} we observe that any metric on $n$
points can be monotonically embedded into an $n$-dimensional Euclidean
space, and that the bound on the dimension is asymptotically tight.
The embedding clearly depends only on the order
of the distances (Lemma \ref{dope-lemma}). We show that for almost every
ordering of the ${n \choose 2}$ distances among $n$
points, the host space of a monotone embedding must be $\Omega(n)$-dimensional.
Similar bounds are given for embeddings into $l_\infty$, and some bounds are also
deduced for other norms.\\
Next we consider embeddings that are even less constrained. Given a metric space
$(X,\delta)$ and some threshold $t$, we seek a mapping $f$ that only respects this
threshold. Namely, $||f(x)-f(y)||<1$ iff $\delta(x,y)<t$.
The input to this problem can thus be thought of as a graph
(adjacency indicating distances below the threshold $t$).
The minimal dimension $d$, such that a graph $G$ can be mapped this way into
$l_2^d$ is known as the {\em sphericity} of G, and denoted $Sph(G)$.
Reiterman, R{\"o}dl and {\v{S}}i{\v{n}}ajov{\'a} (\cite{RRS89a}) show that the
sphericity of $K_{n,n}$ is $n$. This is, then, an explicit example of a metric
space which requires linear dimension to be monotonically embedded into $l_2$.
Other than that, the best lower bounds previously known to us are logarithmic.
In section \ref{prox-graph-sec} we prove a novel lower bound, namely that
for $0 < \delta \leq \half$,
$Sph(G) = \Omega(\frac n {\lambda_2 + 1})$, for any $n$-vertex 
$\delta n$-regular graph, with bounded diameter.
Here $\lambda_2$ is the second largest eigenvalue of the graph.
We also show examples of quasi-random graphs of logarithmic sphericity.
This is somewhat surprising since quasi-random graphs tend to behave like
random graphs, yet the latter have linear sphericity.\\
In our search for further examples of graphs of linear sphericity,
we investigate in section \ref{lambda2-cons} families of graphs whose second
eigenvalue is bounded by a constant (for which the aforementioned lower
bound is linear). We show that such graphs are close to being
complete bipartite,
in the sense that one needs to modify only $o(n^2)$ entries in the adjacency 
matrix to get the latter from
the former. As a corollary, we get that for $0 < \delta < \half$, and $\lambda_2$
there are only finitely many $\delta n$-regular graphs with second eigenvalue
at most $\lambda_2$.

\section{\Dope}\label{dope}
\subsection{Definitions}
Let $X = ([n],\delta)$ be a metric space on $n$ points, such that all pairwise distances are distinct.
Let $||\;||$ be a norm on $\R^d$.
We say that
$\phi:X \rightarrow (\R^d,||\;||)$ is a {\em \dope}
if for every $w,x,y,z \in X$,
$\delta(x,y) < \delta(w,z) \Leftrightarrow ||\phi (x) - \phi (y)|| <
||\phi (w) - \phi (z)||$.

We denote by $d(X,||\;||)$ the minimal $t$ such that there exists a \dope from
$X$ to $(\R^t,||\;||)$. We denote by  $d(n,||\;||) = \max_X d(X,||\;||)$,
the smallest dimension to which every $n$ point metric can be mapped monotonically.

The first thing to notice is that we are actually concerned only with the {\em order}
among the distances between the points in the metric space, and not with the actual
distances.
Let $(X,\delta)$ be a finite
metric space, and let $\rho$ be a linear order on $X \choose 2$. We say that
$\rho$ and $(X,\delta)$ are {\em consistent} if for every
$w,x,y,z \in X$, $\delta(x,y) < \delta(w,z) \Leftrightarrow (x,y) <_\rho (w,z)$.

We start with an easy, but useful observation.
\begin{lemma}\label{dope-lemma}
Let $X$ be a finite set. For every strict order relation $\rho$
on $X \choose 2$, there
exists a distance function $\delta$ on $X$, that is consistent with $\rho$.
\end{lemma}
\proof
Let $\{\epsilon_{ij}\}_{(i,j)\in {X \choose 2}}$ be small, non-negative
numbers, ordered as per $\rho$. Define $\delta(i,j) = 1 + \epsilon_{ij}$.
It is obvious that $\delta$ induces the desired order on the distances of
$X$, and, that if the $\epsilon$'s are small, the triangle inequality holds.
\QED

When we later (Section \ref{l2}) use this observation, we refer to it as
a {\em standard $\epsilon$-construction}, where $\epsilon = \max
\epsilon_{ij}$. It is not hard to see that this metric is Euclidean, that is,
the resulting metric can be isometrically embedded into $l_2$, see
Lemma \ref{momad2} below.

We say that an order relation $\rho$ on $[n] \choose 2$ is
{\em realizable} in $(\R^d,||\;||)$ if there exists a metric space $(X,\delta)$
on $n$ points which is consistent with $\rho$, and a \dope $\phi:X \rightarrow \R^d$.
We say that $\phi$ is a realization of $\rho$.
(In other words, $d(n,||\;||)$ is the minimal $d$ such that
any linear order on $[n] \choose 2$ is realizable in $(\R^d,||\;||)$.)

We denote by $J = J_n$ the $n \times n$ all ones matrix,
and by $PSD_n$ the cone of real symmetric $n \times n$ positive semidefinite
matrices. We omit
the subscript $n$ when it is clear from the context.

Finally, for a graph $G$, and $U, V$ subsets of its vertices, we denote by
$e(U,V) = |\{(u,v) \in E(G): u \in U, v\in V\}|$, and
$e(U) = |\{(u,u') \in E(G): u,u' \in U\}|$.

\subsection{\Dope into $l_\infty$.}
\begin{lemma}
$\frac n 2 - 1 \leq d(n,l_\infty) \leq n$
\end{lemma}
\proof
It is well known that any metric $X$ on $n$ points can be embedded
into $l_\infty^n$ isometrically, hence
$d(n,l_\infty) \leq n$.

For the lower bound, we define
a metric space $(X,\delta)$ with $2n+2$ points
that cannot be realized in $l_\infty^n$.
By lemma \ref{dope-lemma}, it suffices to define an ordering
on the distances. In fact, we define only a partial order, any linear
extension of which will do.
The $2n+2$ points come in $n+1$ pairs, $\{x_i,y_i\}_{i=1,\ldots,n+1}$.
If $z \notin \{x_i,y_i\}$, we let
$\delta(x_i,y_i) > \delta(x_i,z), \delta(y_i,z)$. Assume for contradiction that a
\dope $\phi$ into
$l_\infty^n$ does exist. For each pair $(x,y)$ define
$j(x,y)$ to be some
index $i$ for which $|\phi(x)_i - \phi(y)_i|$ is maximized, that is,
an index $i$ for which  $|\phi(x)_i - \phi(y)_i|=\|\phi(x)-\phi(y)\|_{\infty}$.

By the pigeonhole principle there exist two pairs, say $(x_1,y_1)$ and
$(x_2,y_2)$, for which $j(x_1,y_1) = j(x_2,y_2)=j$.
It is easy to verify that our assumptions on the four real numbers
$\phi(x_1)_j$, $\phi(x_2)_j$,
$\phi(y_1)_j$, $\phi(y_2)_j$, are contradictory.
Thus $d(n,l_\infty) \geq \frac n 2 - 1$.
\QED

\subsection{\Dope into $l_2$.}\label{l2}
\begin{lemma}\label{momad2}
$\frac n 2 \leq d(n,l_2) \leq n$.
Furthermore,  for every  $\delta > 0$, and every large enough $n$, almost no linear orders $\rho$ on ${[n] \choose 2}$ can be realized in dimension less than $\frac n {2+\delta}$.
\end{lemma}
\begin{note}
The upper bound is apparently folklore. As we could not find a reference for it, 
we give a proof here.
\end{note}
\proof
Let $\rho$ be a linear order on ${[n] \choose 2}$. Let $\epsilon$ be a real symmetric
matrix with the following properties:
\begin{itemize}
\item
$\epsilon_{ii} = 0$ for all $i$.
\item
$\frac{1}{n} > \epsilon_{ij} > 0$, for all $i \neq j$.
\item
The numbers $\epsilon_{i,j}$ are consistent with the order $\rho$.
\end{itemize}
Since the sum of each row is strictly less than one, all eigenvalues of $\epsilon$
are in the open interval $(-1,1)$.
It follows that the matrix $I - \epsilon$ is positive definite.
Therefore, there exists a matrix $V$ such that $V V^t = I - \epsilon$.
Denote the $i$'th row of $V$ by $v_i$. Clearly, the $v_i$'s are unit vectors,
and $<v_i,v_j> = - \epsilon_{i,j}$ for $i \neq j$. Therefore,
$||v_i - v_j||_2^2 = <v_i,v_i> + <v_j,v_j> - 2<v_i,v_j> = 2 + 2 \epsilon_{i,j}$.
It follows that the map $\phi(i) = v_i$ is a realization of $\rho$,
and the upper bound is proved. In fact, one can add another point
without increasing the dimension, by mapping it to $0$, and perturbing the
diagonal. \\
For the lower bound, it is essentially known that if $X$ is the metric induced by $K_{n,n}$,
then $d(X,l_2) \geq n$. We discuss this in more detail in the next section.

For the second part of the lemma we need a bound on the number of
{\em sign-patterns} of a sequence of real polynomials. Let
$p_1,...,p_m$ be real polynomials in $l$ variables of (total) degree $d$, and
let $x \in \R^l$ be a point where none of them vanish. The
sign-pattern at $x$ is $(\sgn(p_1(x)),...,\sgn(p_m(x)))$. Denote the
total number of different sign-patterns that can be obtained from $p_1,...,p_m$
by $s(p_1,...,p_m)$. A variation of the Milnor-Thom theorem \cite{Milnor}
due to Alon, Frankl and R\"{o}dl \cite{AFR85} shows:
\begin{theorem}\label{Milnor}\cite{AFR85}
Let $p_1,...,p_m$ be real polynomials as above.
Then for any integer $k$ between 1 and $m$:
\begin{eqnarray*}
s(p_1,...,p_m) \leq 2kd \cdot (4kd - 1) ^ {l + \frac m k - 1}
\end{eqnarray*}
\end{theorem}
Set $n = c \cdot d$, for some constant $c$, and $l = n \cdot d$.
Consider a point $x \in R^l$, and think of it as
an $n \times d$ matrix. Denote the $i$th row of this matrix by $x_i$. As before,
$x$ {\em realizes} an order $\rho$ on ${[n] \choose 2}$ if the distances
$||x_i - x_j||$ are consistent with $\rho$.

For two different pairs, $(i_1,j_1)$ and $(i_2,j_2)$, define the polynomial
\begin{eqnarray*}
p_{(i_1,j_1),(i_2,j_2)}(x) = ||x_{i_1} - x_{j_1}||^2 -
||x_{i_2} - x_{j_2}||^2.
\end{eqnarray*}
The list contains $m = {{n \choose 2} \choose 2}$ polynomials of degree 2.
Note that there is a $1:1$ correspondence between orders on
${[n]} \choose {2}$ and sign-patterns of $p_1,...,p_m$, thus no more
than $s=s(p_1,...,p_m)$ orders may be realized in $l_2^d$.

Take $k=\mu n^2$, for some large constant $\mu$. Then $\log s$ is approximately
$2 c d^2 \log d$. By contrast, that total number of orders is ${n \choose 2}!$,
so its log is about $c^2d^2 \log d$. If $c$ is bigger than 2,
almost all order relations can not be realized.
\QED
\begin{note}
In fact, the same proof shows that for any positive integer $t$,
almost all orders on $n \choose 2$ require linear dimension to be realized,
and in particular that $d(n,l_{2t}) = \Omega(n)$ (where the constant of proportionality 
depends only t):
Simply repeat the argument above with polynomials of degree $2t$
rather than quadratic polynomials. 
\end{note}

\subsection{Other Norms}
We conclude this section
with two easy observations about \dopes into other normed spaces.
The first gives an upper bounds on the dimension required for embedding into $l_p$:
\begin{lemma}
$d(n,l_p) \leq {n \choose 2}$.
\end{lemma}
\proof
By Lemma \ref{momad2}, any metric space on $n$ points can be mapped
monotonically into $l_2$.
It is known (see \cite{DeLa} and also chapter 15 of \cite{Matousek}) 
that any $l_2$ metric on $n$ points can be
isometrically embedded into ${n \choose 2}$-dimensional $l_p$. The composition of these
mappings is a monotone mapping of the metric space into ${n \choose 2}$-dimensional
$l_p$.
\QED

The second observation gives a lower bound for arbitrary norms.
We first note the following:
\begin{lemma}
Let $||\;||$ be an arbitrary $n$-dimensional norm and
let $x_1,...,x_{5^n}$ be points in $\R^n$, such that $||x_i-x_j||>1$
for all $i \neq j$. Then there exits a pair $(x_i,x_j)$ that
$||x_i-x_j|| \geq 2$
\end{lemma}
\proof
Denote by $v$ the volume of $B$, the unit ball in $(\R^n,||\;||)$.
The translates $x_i + \half B$ are obviously
non intersecting, so the volume of their union is $(\frac{5}{2})^n v$.
Assume for contradiction that all pairwise
distances are less than $2$, then all these balls are contained in a single
ball of radius less than $\frac 5 2$. But this is impossible, since the volume of
this ball is less than $(\frac{5}{2})^n v$.
\QED
Note that the $l_\infty$ norm shows that indeed an exponential number of
points is required for the lemma to follow. We do not know, however,
the smallest base of the exponent for which the claim holds. The
determination of this number seems to be of some interest.

\begin{corollary}
There exists an $n$-point metric spaces $(X,\delta)$ such that for any norm $||\;||$,
$d(n,||\;||) = \Omega(\log n)$.
\end{corollary}
\proof
We construct a distance function on $5^n+1$ points
which can not be realized in any $n$-dimensional norm.
By lemma \ref{dope-lemma} it suffices to define a partial order on the distances.
Denote the points in the metric space $0,\ldots,5^n$. Let the distance
between $0$ and any other point be smaller than any distance between any two
points $i \neq j > 0$.
Consider a monotone map $\phi$ of the metric space into $n$-dimensional normed space.
Assume, w.l.o.g., that $\min_{i,j=1,\ldots,5^n}||\phi(i)-\phi(j)|| = 1$.
By the previous lemma there exists a pair of points, $i,j \neq 0$, such that
$||\phi(i)-\phi(j)||>2$. But for $\phi$ to be monotone it must satisfy
$||\phi(0)-\phi(i)||<1$ and $||\phi(0)-\phi(j)||<1$, contradicting the
triangle inequality.
\QED

\section{Sphericity}\label{prox-graph-sec}
So far we have concentrated on embeddings of a metric space into a normed space,
that preserve the order relations between distances. However, in the examples
that gave us the lower bounds for $l_\infty$ and for arbitrary norms,
we actually only needed to distinguish between "long" and "short" distances.
This motivates the introduction of a broader class of maps, that need only
respect the distinction between short and long distances.
More formally, let $X=([n],\delta)$ be a metric space. Its {\em proximity graph} with
respect to some threshold $\tau$, is a graph on $n$ vertices, with an edge between $i$
and $j$ iff $\delta(i,j) \leq \tau$. An embedding of a proximity graph,
is a mapping $\phi$ of its
vertices into normed space, such that $||\phi (i) - \phi (j) || < 1$ iff $(i,j)$
is an edge in the proximity graph (We assume that no distance is exactly 1).
The minimal dimension in which a graph can be so embedded (in Euclidean space)
was first studied by Maehara in \cite{Maehara} under the name {\em sphericity}, and
denoted $Sph(G)$. Following this terminology, we call such an embedding {\em spherical}.\\
The sphericity of graphs  was further studied by Maehara and Frankl in \cite{FraMa},
and then by Reiterman,
R{\"o}dl and {\v{S}}i{\v{n}}ajov{\'a} in \cite{RRS89a}, \cite{RRS89b}, \cite{RRS92}.
Breu and Kirkpatrick have shown in \cite{BK93} that it is NP-hard to recognize
graphs of sphericity 2 (also known as {\em unit disk graphs}) and graphs of sphericity 3.
We refer the reader to \cite{RRS89b} for a
survey of most known results regarding this parameter,
and mention only a few of them here.\\

\begin{theorem}\label{list}
Let $G$ be graph on $n$ vertices with minimal degree $\delta$.
Let $\lambda_n$ the least eigenvalue of its adjacency matrix.
\begin{enumerate}
\item $Sph(K_{m,n}) \leq m + \frac n 2 - 1$ \cite{Maehara}.
\item $Sph(G) = O(\lambda_n^2 \log n)$ \cite{FraMa}.\label{it-lam}
\item $Sph(G) = O((n-\delta) \log (n-\delta))$ \cite{RRS89b}.
\item $Sph(K_{n,n}) \geq n$ \cite{RRS89a}.
\item All but a $\frac 1 n$ fraction of graphs on $n > 37$ vertices have
sphericity at least $\frac n {15} - 1$ \cite{RRS89b}.
\item $Sph(G) \geq \frac {\log \alpha(G)} {\log (2r(G)+1)}$,
where $\alpha(G)$ is the independence number of $G$, and $r(G)$
is its radius \cite{RRS89a}.
\end{enumerate}
\end{theorem}

The first thing to notice is that any lower bower on the sphericity of some
graph on $n$ vertices is also a lower bound on $d(n,l_2)$.
In particular, the fact that $Sph(K_{n,n}) \geq n$ proves the lower bound in
Lemma \ref{momad2}.
(Similarly, any upper bound on the former also applies to the latter.)\\
In this section we are interested in graphs of large sphericity.
The above results tell us that they exist in abundance, yet that graphs
of very small or very large degree have small sphericity
(the maximal degree is an upper bound on $|\lambda_n|$, hence by
(\ref{it-lam}) the sphericity is small if all degrees are small).
Other than the complete bipartite graph, the above results do not point
out an explicit graph with super-logarithmic sphericity.

\subsection{Upper Bound on Margin}
Following Frankl and Maehara \cite{FraMa}, consider an embedding of a
proximity graph where there is a
large margin between short and long distances. In such a situation, the
Johnson-Lindenstrauss Lemma (\cite{JoLi84}) would
yield a spherical embedding into lower dimension: It allows reducing the dimension at the
cost of some distortion. If the distortion is small with respect to the margin, the short and
long distances remain separated.
Alas, we show that for most regular graphs this margin is not large
enough for the method to be useful:
\begin{theorem}
Let $G$ be a $\delta n$-regular graph, with second eigenvalue $\lambda_2 > \frac 2 n$. Let $\phi$ be
an embedding of $G$ as a proximity graph. Denote $a=\max_{u \sim v}||\phi(u) - \phi(v)||_2^2$,
and $b=\min_{u \not \sim v}||\phi(u) - \phi(v)||_2^2$. 
Then $b-a = O(\frac {\lambda_2 + \delta} {\delta n})$.
\end{theorem}
\proof
Denote $m=\min\{1-a, b-1\}$, and for a vertex $i$, denote $v_i = \phi(i)$. The largest value $m$ can attain, over all embeddings $\phi$, is given by the following quadratic semidefinite program:
\begin{eqnarray*}
& \max m & \\
& s.t. \forall (i,j) \in E(G) & ||v_i - v_j||^2 \leq 1-m\\
& \forall (i,j) \notin E(G) & ||v_i - v_j||^2 \geq 1+m
\end{eqnarray*}
Its dual turns out to be:
\begin{eqnarray*}
& \min \half tr A &\\
& s.t. & A \in PSD \\
& \forall (i,j) \in E(G) & A_{ij} \leq 0 \\
& \forall (i,j) \notin E(G), i\neq j & A_{ij} \geq 0 \\
& \forall i & \sum_{j=1,...,n} A_{ij} = 0 \\
& & \sum_{i \neq j} |A_{ij}| = 1
\end{eqnarray*}
Equivalently, we can drop the last constraint, and change the objective function
to $ \min \frac {tr A} {\sum_{i \neq j} |A_{ij}|} $. Next we construct
an explicit feasible solution for the dual program, and conclude from it
a bound on m.

Let $M$ be the adjacency matrix of $G$. Define $A = I + \alpha J - \beta M$.
To satisfy the constraints we need:
\begin{eqnarray*}
& & A \in PSD \\
& & \beta \geq \alpha \geq 0\\
& & 1 + \alpha n - \beta \delta n = 0\\
\end{eqnarray*}
The last condition implies $\alpha = \beta \delta - \frac 1 n$, so it follows
that $\beta \geq \alpha$, and the constraint on $\beta$ is
$\beta \geq \frac {1} {\delta n}$.

Now, since we assume that the graph is $\delta n$-regular, its Perron
eigenvector is $\vec{1}$, corresponding to eigenvalue $\delta n$.
Therefore, we can
consider the eigenvectors of $M$ to be eigenvectors of $J$ and $I$ as well,
and hence also eigenvectors of $A$. If $\lambda \neq \delta n$
is an eigenvalue of $M$, then $1-\beta \lambda$ is an eigenvalue of $A$,
corresponding to the same eigenvector. Denote by $\lambda_2$ the second
largest eigenvalue of $M$, then in order to satisfy the condition $ A \in PSD$ it is enough to set
$\beta = \frac {1} {\lambda_2}$, in which case all the constraints are fulfilled.

We conclude that:
\begin{eqnarray*}
m & \leq & \frac {tr A} {\sum_{i \neq j} |A_{ij}|} =
\frac {n(1+\alpha)}
{\delta n^2 (\beta - \alpha) + ((1-\delta) n^2 - n) \alpha} \\
& = &
\frac {n+ \frac {\delta n} {\lambda_2} - 1}
{\delta n (\frac {n+\delta n} {\lambda_2} - 1) + ((1-\delta) n - 1)
(\frac {\delta n} {\lambda_2} - 1)}
 <
4 \frac {1 + \frac {\delta} {\lambda_2}} {\frac {\delta n} {\lambda_2}}
 = 4 \frac {\lambda_2 + \delta} {\delta n}.
\end{eqnarray*}
In particular, $b-a = O(\frac {\lambda_2 + \delta} {\delta n})$.
\QED

In order to derive a non trivial result from Johnson-Lindenstrauss lemma,
we need that $\frac 1 {m^2} \log n = o(n)$, and in particular that
$\lambda_2=\Omega(\delta \sqrt{n \log n})$. The above shows that this can
happen only if $\lambda_2 = \omega(\delta \sqrt{n \log n})$.
On the other hand, Frankl and Maehara show that their method does give a non trivial
bound when $\lambda_n = o(\sqrt{\frac n {\log n}})$.
Consequently, we get that a $\delta n$-regular graph (think of $\delta$ as constant)
can't have both $\lambda_2=o(\sqrt{n \log n})$ and
$\lambda_n = o(\sqrt{\frac n {\log n}})$. This is a bit more subtle than what
one gets from the second moment argument, namely, that the graph can't have both
$\lambda_2=o(\sqrt{n})$ and $\lambda_n = o(\sqrt{n})$.

\subsection{Lower Bound on Sphericity}

\begin{theorem}\label{our-bound}
Let $G$ be a $d$-regular
graph with diameter $D$ and $\lambda_2$,
the second largest
eigenvalue of $G$'s adjacency matrix, is at least $d - \half n$.
Then $Sph(G) = \Omega(\frac {d - \lambda_2} {D^2(\lambda_2 + O(1))})$.
\end{theorem}

In the interesting range where $d \leq \frac n 2$,
and $\lambda_2 \geq 1$
the bound is
$Sph(G) = \Omega(\frac {d - \lambda_2} {D^2 \lambda_2})$.
\proof
It will be useful to consider the following operation on matrices. Let $A$
be an $n \times n$ symmetric matrix, and denote by $\arrow{a}$ the vector
whose $i$-th coordinate
is $A_{ii}$. Define $R(A)$ to be the $n \times n$ matrix with all
rows equal to $\arrow{a}$, and $C(A) = R(A)^t$. Define:
\begin{eqnarray*}
\breve{A} = 2A - C(A) - R(A) + J
\end{eqnarray*}
First note that the rank of $\breve{A}$ and that of $A$ can differ by
 at most 3.
%Now, suppose that $A$ is positive semi-definite,
%then the restriction of the linear map $\breve{A}$ to the subspace $\vec{1} ^ \bot$
%is positive semi-definite,
%since for any $x$ in this subspace, $Jx = xR(A) = C(A)x = 0$.
Now, consider the case where $A$ is the Gram matrix of some vectors
 $v_1,...,v_n \in R^d$.
Then all diagonal entries of $\breve{A}$ equal one, and the $(i,j)$ entry is
2$<v_i,v_j> - <v_i,v_i> - <v_j,v_j> + 1 = 1 - ||v_i-v_j||^2$.

We will need the following lemma (see \cite{HoJo}, p.175):
\begin{lemma}\label{d-lemma}
Let $X$ be a real symmetric matrix, then
$rank(X) \geq \frac {(tr X)^2} {\sum_{i,j}X_{i,j}^2}$
\end{lemma}

Applying this to $\breve{A}$, we conclude that:
\begin{eqnarray}\label{d-bound}
rank(\breve{A}) \geq \frac {n^2} {n +
\sum_{i \neq j} (1 - ||v_i - v_j||^2)^2 }
\end{eqnarray}

Let $v_1,...,v_n \in \R^d$ be an embedding of $G$. By the discussion above
it is enough to show that
\begin{eqnarray}\label{denom}
\sum_{i \neq j} (1-||v_i - v_j||^2)^2 = O(D^2 n^2 \frac {\lambda_2}
{d - \lambda_2}).
\end{eqnarray}
By the triangle inequality $||v_i - v_j|| \leq D$ for any two vertices.
So the LHS of (\ref{denom}) is bigger by at most a factor of $D^2$ than:
\begin{eqnarray*}
&& \sum_{(i,j) \notin E} (||v_i - v_j||^2 - 1) +
\sum_{(i,j) \in E} (1 - ||v_i - v_j||^2) =
\end{eqnarray*}
\begin{eqnarray}\label{hs-norm}
&& \sum_{(i,j) \notin E} ||v_i - v_j||^2  -
\sum_{(i,j) \in E} ||v_i - v_j||^2 - {n \choose 2} + nd
\end{eqnarray}
We can bound this sum from above, by solving the following SDP:
\begin{eqnarray*}
& \max & \sum_{(i,j) \notin E} (V_{ii} + V_{jj} - 2V_{ij}) +
\sum_{(i,j) \in E} (- V_{ii} - V_{jj} + 2V_{ij}) - {n \choose 2} + nd\\
& s.t. & V \in PSD \\
& \forall (i,j) \in E & V_{ii} + V_{jj} - 2V_{ij} \leq 1 \\
& \forall (i,j) \notin E & V_{ii} + V_{jj} - 2V_{ij} \geq 1
\end{eqnarray*}
The dual problem is:
\begin{eqnarray*}
& \min & \half tr A \\
& s.t. & A \in PSD \\
& \forall (i,j) \in E & A_{ij} \leq -1 \\
& \forall (i,j) \notin E, i\neq j &  A_{ij} \geq 1 \\
& \forall i \in [n] & \sum_{j=1,...,n} A_{ij} = 0
\end{eqnarray*}
Let $M$ by the adjacency matrix of the graph, and set
$A = (\alpha d - n)I + J - \alpha M$, where $\alpha \geq 2$
will be determined shortly. This takes care of the all constraints
except for $A \in PSD$. Note that since $M$ is regular,
its eigenvectors
are also eigenvectors of $A$.
Moreover, if $M u = \lambda u$ for a non constant $u$, then
$A u =  \alpha d - n - \alpha \lambda$ (and $A \vec{1} = 0$).
% and that for any such vector other than
%$\vec{1}$, if in $M$ it is associated with eigenvalue $\lambda$, in $A$ it
%corresponds to $\alpha d - n - \alpha \lambda$ ($\vec{1}$ corresponds to 0).
So take
$\alpha = \frac {n} {d - \lambda_2}$, and by our assumption on $\lambda_2$,
$\alpha \geq 2$.

Now $A$ gives an upper bound on (\ref{hs-norm}):
\begin{eqnarray*}
\half tr A =
\half n(\alpha d - n + 1) = \half n^2 \frac d {d-\lambda_2} - \half n^2 
+ \half n= \half n^2 \frac {\lambda_2} {d-\lambda_2} + \half n.
\end{eqnarray*}
This, by (\ref{d-bound}), shows that the dimension of the embedding is
$\Omega\left(\frac {d - \lambda_2} {D^2(\lambda_2 + O(1))} \right)$.
\QED

\subsection{A Quasi-random Graph of logarithmic Sphericity}
It is an intriguing problem to construct new examples of graphs of linear
sphericity. Since random graphs have this property, it is natural to search
among quasi-random graphs.
There are several equivalent definitions for such graphs (see \cite{AlSp}). The one we adopt here is:
\begin{definition}
A family of graphs is called {\em quasi-random} if the graphs in the family are
 $(1+o(1))\frac n 2$-regular, and all their
eigenvalues except the largest one are (in absolute value) $o(n)$ .
\end{definition}

Counter-intuitively, perhaps, quasi-random graphs may have very small sphericity.
\begin{lemma}\label{alex}
Let $\G$ be the family of graphs with vertex set $\{0,1\}^k$, and edges connecting vertices that are at Hamming distance at most $\frac k 2$.
Then $\G$ is a family of quasi-random graphs of logarithmic sphericity.
\end{lemma}
\proof
The fact that the sphericity is logarithmic is obvious -
simply map each vertex to
the vector in $\{0,1\}^n$ associated with it. To show that all eigenvalues except
the largest one are $o(2^k)$ we need the following facts about Krawtchouk
polynomials (see \cite{vL}). Denote by
$K_s^{(k)}(i) = \sum_{j=0}^s(-1)^j{i \choose j}{{k-i} \choose {s-j}}$ the
Krawtchouk polynomial of order $s$ over $\Z_2^k$.
For simplicity we assume that $k$ is odd.
\begin{enumerate}
\item For any $x \in \Z_2^k$ with $|x|=i$,
$\sum_{z \in \Z_2^k\\|z|=s}(-1)^{<x,z>} = K_s^{(k)}(i)$.
\item $\sum_{s=0}^l K_s^{(k)}(i) = K_l^{(k-1)}(i-1)$.
\item For any $s$ and $k$, $\max_{i=0,\dots,n} |K_s^{(k)}(i)| = K_s^{(k)}(0) = {k \choose s}$.
\end{enumerate}
Observe that $G$ is a Cayely graph for the group $\Z_2^k$ with generator
set $\{g \in \Z_2^k : |g| \leq \frac k 2\}$. Since $\Z_2^k$ is abelian, the
eigenvectors of the graphs are independent of the generators, and are simply the
characters of the group written as the vector of their values. Namely,
corresponding each
$y \in \Z_2^k$ we have an eigenvector $v^y$, such that $v^y_x = (-1)^{<x,y>}$.
For every $y$, $v^y_0 = 1$, so to figure out the eigenvalue
corresponding to $v^y$, we simply need to sum the value of $v^y$ on the
neighbors of $0$. Note that for $y=0$ we get the all $1$s vector, which corresponds to the largest eigenvalue. So we are interested in $y$'s such that $|y| > 0$. By the first two facts above we have:
\[\lambda_y = \sum_{g \in \Z_2^k, |g| \leq \frac k 2}(-1)^{<y,g>} = \sum_{s=0}^{\frac {k-1} 2} K_s^{(k)}(|y|) = K_{\frac {k-1} 2}^{(k-1)}(|y|-1).\]
By the third fact, this is at most ${{k-1} \choose {\frac {k-1} 2}}
\approx \frac {2^{k-1}} {\sqrt{k-1}} =
o(2^{k-1})$.
\QED

\section{Graphs with bounded $\lambda_2$}\label{lambda2-cons}
Theorem \ref{our-bound} suggests families of graphs that have linear sphericity.
Namely, for $0 < \delta \leq \half$, and $\lambda_2 > 0$, the theorem says that
$\delta n$-regular graphs with second eigenvalue at most $\lambda_2$ have
linear sphericity. In this section we characterize such graphs. We prove
that for $\delta = \half$ such graphs are nearly complete bipartite, and that for
other values, only finitely many graphs exist.\\
It is worth noting that graphs with bounded second eigenvalue have been 
previously studied. The apex of these works is probably that of
Cameron, Goethals, Seidel and Shult, who characterize in \cite{CGSS} graphs with
second eigenvalue at most 2.

\subsection{$n/2$-regular graphs} 
In this section we consider the family $\G$ of $n/2$-regular graphs,
and second largest eigenvalue $\lambda_2$ bounded by a constant.
We prove that, asymptotically, they are nearly complete bipartite.

\begin{definition}
Let $G$ and $H$ be two graphs on $n$ vertices. We say that $G$ and $H$ are {\em close}, if there is a labeling of their vertices such that $|E(G) \bigtriangleup E(H)| = o(n^2)$.
\end{definition}

\begin{theorem}\label{main1}
Every $G \in \G$ is close to $K_{n/2,n/2}$, where $n$ is the number of vertices in $G$.
\end{theorem}
\begin{note}
By passing to the complement graph, if $\lambda_n = O(1)$,
then $G$ is close to the disjoint union of two cliques, $K_{n/2} \dot{\cup} K_{n/2}$.
\end{note}

We need several lemmas. The first is the well-known expander mixing lemma \cite{FriPi}. The second
is a special case of Simonovitz's stability theorem (\cite{Sim}),
for which we give a simple proof here. The third is a commonly used corollary of
Szemeredi's Regularity Lemma. We shall also make use of the Regularity Lemma itself
(see e.g. \cite{Diestal}).

\begin{lemma}\label{no-clique}
Let $G$ be an $\frac n 2$-regular graph on $n$ vertices with second largest
eigenvalue $\lambda_2$. Then every subset of vertices with $k$ vertices has at
most $\frac 1 4 k^2 + \half \lambda_2 k$ internal edges.
\end{lemma}

\begin{lemma}\label{close-bi}
Let $R$ be a triangle-free graph on $n$ vertices, with $n^2/4 - o(n^2)$ edges.
Then $R$ is close to $K_{n/2,n/2}$. Furthermore, all but $o(n)$ of the vertices
have degree $\frac n 2 \pm o(n)$.
\end{lemma}
\proof
Denote by $d_i$ the degree of the $i$th vertex in $R$, and by $m$ the number of edges.
Then:
\[\sum_{(i,j) \in E(R)}(d_i+d_j) = \sum_{i \in V(R)} d_i^2 \geq \frac 1 n
(\sum_{i \in V(R)} d_i)^2 = \frac {4m^2} n.\]
Thus, there is some edge $(i,j) \in E(R)$ such that
$d_i + d_j \geq \frac {4m} n = n - o(n)$. Let $\Gamma_i$ and $\Gamma_j$ be the
neighbor sets of $i$ and $j$. Since $i$ and $j$ are adjacent, and $R$ has no triangles,
the sets $\Gamma_i$ and $\Gamma_j$ are disjoint and independent.
If we delete the $o(n)$ of vertices in $V\backslash (\Gamma_i \cup \Gamma_j)$
we obtain a bipartite graph. We have deleted only $o(n^2)$ edges, so the remaining
graph still has $n^2/4 - o(n^2)$ edges. But this means that
$|\Gamma_i|, |\Gamma_j| = \frac n 2 - o(n)$, and that the degree of each vertex
in these sets is $\frac n 2 \pm o(n)$
\QED

Recall that the Regularity Lemma states that for every $\epsilon > 0$ and
$m \in \N$ there's an $M$, such that
the vertex set of every large enough graph can be partitioned into $k$ subsets,
for some $m \leq k \leq M$ with the following properties: All subsets except
one, the ``exceptional'' subset, are of the same size. The exceptional subset
contains less than an
$\epsilon$-fraction of the vertices. All but an $\epsilon$-fraction of the pairs of
subsets are $\epsilon$-regular.\\
The regularity graph with respect to such a partition and a threshold $d$,
has the $k$ subsets as vertices. Two subsets, $U_1$ and $U_2$ are adjacent,
if they are $\epsilon$-regular, and $e(U_1,U_2) > dn^2$.

\begin{lemma} [\cite{Diestal}, Lemma 7.3.2]
Let $G$ be a graph on $n$ vertices, $d \in (0,1]$, $\epsilon = d^{-4}$.
Let $R$ be an $\epsilon$-regularity graph of $G$, with (non exceptional) sets of
size at least $\frac s {\epsilon}$, and threshold $d$. If $R$ contains a triangle,
then $G$ contains a complete tripartite subgraph, with each side of size $s$.
\end{lemma}
\begin{corollary}
If $G \in \G$, and $R$ is as in the lemma, with $s = 10 \lambda_2$, then $R$ is
triangle free. In this case, if $R$ has $\frac {k^2} 4 - o(k^2)$ edges,
then $R$ is close to complete bipartite.
\end{corollary}
\proof
If $R$ contains a triangle, then $G$ contains a complete tripartite subgraph,
with $s$ vertices on each side. Let $U$ be the set vertices in this subgraph.
Then $e(U) = 3s^2 = 300 \lambda_2^2$, but by lemma \ref{no-clique}
$e(U) \leq 50 \lambda_2^2$ - a contradiction. The second part now follows from
Lemma \ref{close-bi}.
\QED

\proof (Theorem)
We would like to
apply the Regularity Lemma to graphs in $\G$, and have $\epsilon = o(1)$,
and $k = \omega(1)$ as well as $k = o(n)$. Indeed, this can be done.
Since $M$ depends only on $m$ and $\epsilon$, choose
$d=o(1)$, and $m = \omega(1)$, such that the $M$ given by the lemma satisfies $\frac n {(M+1)} \geq \frac s {\epsilon}$. As $M$ depends only on $m$ and $\epsilon$, $\frac M {\epsilon}$ can be made small enough, even with the requirements
on $d$ and $m$.

Let $R$ be the regularity graph for the partition given by the Regularity Lemma, with threshold $d$ as above. Denote by $k$ the number of sets in the partition, and their size by $l$ (so $k\cdot l = n(1-\eta)$, for some $\eta \leq \epsilon$). We shall show that $R$ is close to complete bipartite, and that $G$ is close the graph obtained by replacing each vertex in $R$ with $l$ vertices, and replacing each edge in $R$ by a $K_{l,l}$.

Call an edge in $G$ $(i)$ ``irregular'' if it belongs to an irregular pair; $(ii)$ ``internal'' if it connects two vertices within the same part; $(iii)$ ``redundant'' if it belongs to a pair of edge density smaller than $d$, or touches a vertex in the exceptional set.  Otherwise $(iv)$, call it ``good''.\\
Recall that $\epsilon = o(1)$, so only $o(k^2)$ pairs of sets are not $\epsilon$-regular. Thus, $G$ can have only $o(l^2k^2) = o(n^2)$ irregular edges. Also, $d = o(1)$, so the number of redundant edges is $k^2 \cdot o(l^2) + o(l) \frac n 2 = o(n^2)$. Finally, the number of internal edges is at most $\half l^2 k$, hence there are $\frac {n^2} 4 - o(n^2)$ good edges. \\
The number of edges between two sets is at most $l^2$, so $R$ must have at least
\[\frac {n^2 - o(n^2)} {4l^2} = \frac {k^2} 4 - o(k^2)\]
edges. The corollary implies that it is close to complete bipartite.
By lemma \ref{close-bi}, the valency of all but $o(k)$ of the vertices in $R$ is indeed $\frac k 2 \pm o(k)$. This means
that every edge in $R$ corresponds to $l^2 - o(l^2)$ good edges in $G$ (as the number of edges in $R$ is also no more than $\frac {k^2} 4 + o(k^2)$).\\
To see that $G$ is close to complete bipartite, let's count how many edges need to be modified. First, delete $o(n^2)$ edges that are not ``good''. Next add all possible $o(n^2)$ new edges between pairs of sets that have ``good'' edges between them. As $R$ is close to complete bipartite, we need to delete or add all edges between $o(k^2)$ pairs. Each such step modifies $l^2$ edges, altogether $o(l^2k^2) = o(n^2)$ modifications. Finally, divide the $o(n)$ vertices of the exceptional set evenly between the two sides of the bipartite graph, and add all the required edges, and the tally remains $o(n^2)$.
\QED

\begin{note}
In essence, the proof shows that a graph with no dense induced subgraphs is close to complete bipartite. This claim is similar in flavor to Bruce Reed's {\em Mangoes and Blueberries} theorem \cite{Reed99}. Namely, that if every induced subgraph $G'$ of $G$ has an independent set of size $\half |G'| - O(1)$, then $G$ is close to being  bipartite. The conclusion in Reed's theorem is stronger in that only a {\em linear} number of edges need to be deleted to get a bipartite graph.
\end{note}
\begin{note}
In fact, the proof gives something a bit stronger. Let $t_r(n)$ be the number
of edges in an $n$-vertex complete $r$-partite graph, with parts of equal size.
Using the general Stability Theorem (\cite{Sim})
instead of Lemma \ref{close-bi}, the same proof shows that if a graph has
$t_n - o(n^2)$ edges and no dense induced subgraphs, then it is close to being
complete $r$-partite.
\end{note}

\subsection{$\delta n$-regular graphs}
In Theorem \ref{main1} we required that the degree is $n/2$. We can deduce from the theorem that
this requirement can be relaxed:
\begin{corollary}
Let $\G$ be a family of $d$-regular graphs, with $d \leq \frac n 2$, ($n$ being
the number of vertices in the graph) and bounded second eigenvalue, then every $G \in \G$ is close to a complete bipartite graph.
\end{corollary}
\proof
Let $M \in \M_n$ be the adjacency matrix of such a $d$-regular graph, and denote $\bar{M} = J - M$, where $J$ is the all ones matrix. Consider the graph $H$ corresponding to the following matrix:
\begin{eqnarray*}
N =
\left(
\begin{array}{cc}
M & \bar{M} \\
\bar{M}^t & M
\end{array}
\right)
\end{eqnarray*}
Clearly $H$ is an $n$-regular graph on $2n$ vertices.
Denote by $(x, y)$ the concatenation of two $n$-dimensional vectors, $x$, $y$, into a $2n$ dimensional vector. Let $v$ be an eigenvector of $M$ corresponding to eigenvalue $\lambda$.
It is easy to see that $v$ is also an eigenvalue of $\bar{M}$: If $v = \vec{1}$ (and thus $\lambda = d$) it corresponds to eigenvalue $n - \lambda$, otherwise to $(-\lambda)$.\\
Thus, $(v, v)$ and $(v, -v)$ are both eigenvectors of $N$. If $v = \vec{1}$ they correspond to eigenvalues $n$, $2d-n$, respectively, otherwise to $0$, $2\lambda$. 
Since the $v$'s are linearly independent, so are 
the $2n$ vectors of the form $(v, v)$ and $(v, -v)$:
Consider a linear combination of these vectors that gives $0$. Both the sum and the difference of the coefficients of each pair have to be $0$, and thus both are 0. So we know the entire spectrum of $N$, and see, since $d \leq \frac n 2$, that theorem \ref{main1} holds for it.\\
Let $H'$ be a complete bipartite graph that is close to $H$. Since $H$ differs from $H'$ by $o(n^2)$ edges, the same holds for subgraphs over the same set of vertices. In particular, $G$ is close to the subgraph of $H'$ spanned by the first $n$ vertices. Obviously, every such subgraph is itself complete bipartite.
\QED

\begin{corollary}\label{no-graphs}
For every $0 < \delta < \half$ and $c$, there are only finitely many $\delta n$-regular graphs with $\lambda_2 < c$.
\end{corollary}
\proof
Consider such a graph with $n$ large. By the previous corollary it is close
to complete bipartite. Since it is also regular, it must be close to
$K_{\frac n 2, \frac n 2}$, which contradicts the constraint $\delta < \half$.
\QED

\subsection{Graphs with both $\lambda_2$ and $\lambda_{n-1}$ bounded by a 
constant}
Theorem \ref{main1} can loosely be stated as follows: A regular graph with 
spectrum similar to that of a bipartite graph ($\lambda_1$ being close to $n/2$
and $\lambda_2$ being close to $0$) is close to being complete bipartite.
We conclude this section by noting that if we strengthen the assumption on how
close the spectrum of a graph is to that of a bipartite graph, we get a stronger
result as to how close it is to a complete bipartite graph. 

\begin{theorem}\label{main2}
Let $\G$ be a family of $\frac n 2$-regular graphs on $n$ vertices, with 
both $\lambda_2$ and $\lambda_{n-1}$ bounded by a 
constant. Then every
$G \in \G$ is close to a $K_{\frac n 2,\frac n 2 }$, in the sense that such a graph can be obtained from $G$ by modifying a linear number of edges for $O(\sqrt{n})$ vertices of $G$, and $O(\sqrt{n})$ edges for the rest.
\end{theorem}
\proof
First note that it follows that $\lambda_n(G) = -\frac n 2 + O(1)$. Take 
$G \in \G$, and let $A$ be its adjacency matrix.
Clearly $tr(A^2) = \frac {n^2} 2$. If $\lambda_{n-1}(G) = - O(1)$, then
\[\frac {n^2} 2 = tr(A^2) = 
\lambda_1^2 + \lambda_n^2 + \sum_{i=2,...,n-1}\lambda_i^2\] Since $\lambda_1 = \frac n 2$
\[\lambda_n^2 = 
\frac {n^2} 2 - \left(\frac n 2\right)^2 - \sum_{i=2,...,n-1}\lambda_i^2\]
As $\lambda_2,\ldots,\lambda_{n-1} = O(1)$ we have
\[\lambda_n^2 = \frac {n^2} 4 + O(n)\]
And since $\lambda_n$ is negative, and is smaller than $\lambda_1$ in absolute value:
\[\lambda_n = -\frac n 2 + O(1).\]

Let $x$ be an eigenvector corresponding to $\lambda_n$. Suppose, w.l.o.g that $||x||_{\infty} = 1$ and that $x_v = 1$. Denote $A = \{u : x_u \leq -(1-\frac 1 {\sqrt{n}})\}$, and $B = \{w : x_w \geq (1-\frac 1 {\sqrt{n}})\}$. The eigenvalue
condition on $v$ entails:
\begin{eqnarray*}
\frac n 2 - O(1) = -\sum_{u:(u,v) \in E} x_u.
\end{eqnarray*}
Thus, there is a vertex $u$ such that $x_u \leq -(1-O(\frac 1 n))$. It is not hard to verify that $v$ must have $\frac n 2 - O(\sqrt{n})$ neighbors in $A$,
and that $u$ must have $\frac n 2 - O(\sqrt{n})$ neighbors in $B$.

Now denote $A' = \{u : x_u \leq -\half\}$, and $B' = \{w : x_w \geq \half \}$.
Again, it is not hard to check that each vertex in $A$ must have
$\frac n 2 - O(\sqrt{n})$ neighbors in $B'$, and vice versa. Thus, delete the
$O(\sqrt{n})$ vertices that are neither in $A$ nor in $B$. For each remaining
vertex in $A$ (similarly in $B$), its degree is at most $\frac n 2$, and
at least $\frac n 2 - O(\sqrt{n})$. It has $\frac n 2 - O(\sqrt{n})$ neighbors
in $B$, so the number of its neighbors in $A$, and the number of its non-neighbors in $B$ is $O(\sqrt{n})$. By deleting and adding $O(\sqrt{n})$ edges to each
vertex, we get a complete bipartite graph.
\QED

\begin{note}
Alternatively, we could have defined $\G$ as a family of $\frac n 2$-regular
graphs with 
$\lambda_2$ bounded, and $\lambda_n(G) = -\frac n 2 + O(1)$. It's interesting
to note that in this case it follows that $\lambda_{n-1}$ is bounded.
For $G \in G$, if $G$ is bipartite, then it is complete bipartite, and $\lambda_{n-1}(G) = 0$. Otherwise, $\chi(G) > 2$, and by a theorem of Hoffman (\cite{Hoff}) $\lambda_n(G) + \lambda_{n-1}(G) + \lambda_1(G) \geq 0$. By our assumption,  $\lambda_n(G) + \lambda_1(G) = O(1)$, and since $\lambda_{n-1}(G) < 0$ (otherwise the eigenvalues won't sum up to 0), it follows that $\lambda_{n-1}(G) = -O(1)$.
\end{note}

\section{Conclusion and Open Problems}

The only explicit examples known so far for graphs that have linear sphericity are $K_{n,n}$ and small modifications of it.
We conjecture that more complicated graphs, such as the Paley graph, also have
linear sphericity. Note that the lower bound presented here only shows a bound
of $\Omega(\sqrt{n})$. It is also interesting to know if the bound can be improved,
either as a pure spectral bound, or with some further assumptions on the
structure of the graph. \\
What is the largest sphericity, $d=d(n)$, of an $n$-vertex graph? We know that
$\frac n 2 \leq d \leq n-1$. Can this gap be closed? For a seemingly related question,
the smallest dimension required
to realize a sign matrix (see \cite{AFR85}) the answer is known to be
$\frac n 2 \pm o(n)$.
We have also seen a similar gap for $d(n,l_2)$ and $d(n,l_\infty)$.
Can this be
closed? Can some kind of interpolation arguments generalize the bounds we know for
these two numbers to bounds on $d(n,l_p)$ for $p>2$?\\
Our interest in sphericity arose from a search for a lower bound on $d(n,l_2)$.
But why limit the discussion to Euclidean space? What can be said of spherical
embeddings into $l_1$ or $l_{\infty}$? The former may be particularly
interesting, as it will give a non-trivial lower bound on $d(n,l_1)$.\\
We have seen that $\frac n 2$-regular graphs with bounded
second eigenvalue are $o(n^2)$-close to complete bipartite. However, the only
example we know of such graphs are constructed by taking a complete bipartite
graph, and changing a constant number of edges for each vertex. These graphs
are $O(n)$-close to being complete bipartite. Are there examples of such
families which are further from complete bipartite graphs, or can a stonger
notion of closeness be proved? 

\section{Acknowledgments}
We would like to thank Alex Samorodnitsky for showing us how to prove Lemma \ref{alex} and Noga Alon for sending us Lemma \ref{d-lemma}.

%%%%%%%%%%%%%%%%%%%%%%%%%%%%%%%%%%%%%%%%%%%%%%%%%%%%%%%%%%%%%%%%%%%

\bibliography{bib}

\end{document}